\documentclass[
11pt,%
tightenlines,%
twoside,%
onecolumn,%
nofloats,%
nobibnotes,%
nofootinbib,%
superscriptaddress,%
noshowpacs,%
centertags]%
{revtex4}

\usepackage{amssymb}
\usepackage{ljm}

\newcommand{\Cc}{\mathbb{C}}
\newcommand{\geqs}{\geqslant}
\newcommand{\leqs}{\leqslant}
\newcommand{\Rr}{\mathbb{R}}

\newtheorem{remark}{Remark}

\begin{document}

\titlerunning{Exponential and backward Euler for nonlinear heat equations}

\authorrunning{Botchev, Zhukov}

\title{Exponential Euler and backward Euler methods for nonlinear heat conduction problems}

\author{\firstname{Mikhail A.}~\surname{Botchev}}
\email[E-mail: ]{botchev@kiam.ru}
\affiliation{Keldysh Institute of Applied Mathematics of Russian Academy of Sciences, 
Miusskaya Sq., 4, 125047 Moscow, Russia} 
\affiliation{Marchuk Institute of Numerical Mathematics of Russian Academy of Sciences, 
Gubkina St., 8, 119333 Moscow, Russia} 

\author{\firstname{Victor T.}~\surname{Zhukov}}
\email[E-mail: ]{vic.zhukov@gmail.com}
\affiliation{Keldysh Institute of Applied Mathematics of Russian Academy of Sciences, 
Miusskaya Sq., 4, 125047 Moscow, Russia} 

\firstcollaboration{(Submitted by \dots)}


\received{September 1, 2022}

\begin{abstract}
In this paper a variant of nonlinear exponential Euler scheme
is proposed for solving  nonlinear heat conduction problems.
The method is based on nonlinear iterations where at each iteration 
a linear initial-value problem has to be solved.  
We compare this method to the backward Euler method
combined with nonlinear iterations.  
For both methods we show 
monotonicity and boundedness of the solutions and give sufficient 
conditions for convergence of the nonlinear iterations.
Numerical tests are presented to examine performance of 
the two schemes.
The presented exponential Euler scheme is implemented
based on restarted Krylov subspace methods and, hence, is 
essentially explicit (involves only matrix-vector products). 
\end{abstract}
\maketitle

\section{Introduction and problem setting}
We consider a class of nonlinear heat conduction problems
\begin{equation}
\label{heat}
\frac{\partial}{\partial t}u(\bm{x},t) = \nabla\cdot (k(u)\nabla u(\bm{x},t))
+g(\bm{x},t),
\quad
u(\bm{x},0) = u^0(\bm{x}),
\quad 
u(\bm{x},t) \bigr|_{\partial\Omega}=b(\bm{x},t),
\end{equation}
where $\bm{x}\in\Omega\subset\Rr^d$, $d\in\{1,2,3\}$, $t\in[0,T]$,
$\nabla u$ denotes the gradient of $u(\bm{x},t)$, $\nabla\cdot$ is 
the divergence operator and functions $u^0(\bm{x})$, $b(\bm{x},t)$ and 
$g(\bm{x},t)$ are given.
Furthermore, we have 
\begin{equation}
\label{k}
k(u) = k_0 u^\sigma,  
\end{equation}
with given constants $k_0>0$ and $\sigma>0$, and it holds
\begin{equation}
\label{pos}
b(\bm{x},t)\geqs 0, \quad 
g(\bm{x},t)\geqs 0, \quad
u^0(\bm{x})\geqs 0, \quad 
\bm{x}\in\Omega,\; t\geqs 0.
\end{equation}
We assume that a suitable spatial discretization of initial-boundary problem (IVP)
in~\eqref{heat} leads to initial-value problem 
\begin{equation}
\label{ivp}
y'(t) = -A(y(t)) y(t) + g(t),\quad y(0)=v, \quad\text{with}\; v\geqs 0,
\end{equation}
where the matrix $A(y)\in\Rr^{N\times N}$ is symmetric positive semidefinite for
any $y\in\Rr^N$ and its off-diagonal entries are nonpositive, i.e., 
\begin{equation}
\label{aij}
a_{ij}(y)\leqs 0, \quad i\ne j, \quad \forall y\in\Rr^N.  
\end{equation}
In~\eqref{ivp} and throughout the paper we understand vector inequalities elementwise,
i.e., $v\geqs 0$ means that all the entries of the vector $v\in\Rr^N$ are 
nonnegative. 
Note that the function $g: \Rr\rightarrow\Rr^N$ in~\eqref{ivp} contains not only the values 
of the source function $g(\bm{x},t)$ on the grid but also the grid contributions from 
the boundary conditions.  These contributions will typically be of the form
$-a_{ij}(y)b(\bm{x},t)$ for some $i\ne j$ and $\bm{x}\in\partial\Omega$.
In view of~\eqref{pos},\eqref{aij}, we may assume that the values of the function $g(t)$ 
are nonnegative vectors in $\Rr^N$, i.e.,
\begin{equation}
\label{g>=0}
g(t)\geqs 0, \quad t\geqs 0. 
\end{equation}
Furthermore, throughout the paper $(x,y)=y^Tx$, $x,y\in\Rr^N$, denotes the standard inner 
product and, unless reported otherwise, $\|\cdot\|$ is the Euclidean vector or matrix norm.
We assume that IVP~\eqref{ivp} has a unique solution and
there is a constant $L>0$  
\begin{equation}
\label{L}
\|A(u) - A(v)\| \leqs L\|u-v\|, \quad \forall u,v\in\Rr^N.
\end{equation}

\begin{corollary}
\label{IVPpos}
Solution $y(t)$ of the semidiscrete IVP~\eqref{ivp} is entrywise nonnegative, i.e.,
$$
y(t)\geqs 0, \quad t\geqs 0 
$$
provided that conditions~\eqref{aij} and~\eqref{g>=0} hold.  
\end{corollary}

\begin{proof}
The statement of the corollary follows from considering the explicit Euler scheme
applied to~\eqref{ivp}:
$$
y^{n+1} = y^n - \Delta t A (y^n) y^n + g^n, 
$$  
where the superscript $\cdot^n$ denotes the time step number and
$\Delta t>0$ is the time step size.  The last relation reads elementwise,
for each entry $y_i^{n+1}$, $i=1,\dots, N$, of the vector $y^{n+1}$ as
$$
y_i^{n+1} = (1 - \Delta t a_{ii}(y^n))y_i^n - \sum_{j\ne i} a_{ij}(y^n) y_j^n + g_i^n.
$$ 
Since $1 - \Delta t a_{ii}(y^n)\geqs 0$ for sufficiently small $\Delta t>0$ and 
due to~\eqref{aij} and~\eqref{g>=0}, we have $y_i^{n+1}\geqs 0$, $n=1,2,\dots$.
For $\Delta t\rightarrow 0$ this nonnegative numerical solution converges to
the unique solution of~\eqref{ivp}.  
\end{proof}

\section{Solution methods}
\subsection{The backward Euler method}
The backward Euler method applied to~\eqref{ivp} reads 
\begin{align}
\frac{y^{n+1}-y^n}{\Delta t} &= -A(y^{n+1}) y^{n+1} + g^{n+1} \quad\Leftrightarrow\quad
\notag
\\
\label{BE}
(I+\Delta t A(y^{n+1}))y^{n+1} &= y^n + \Delta t g^{n+1},
\quad n=0,1,2,\dots.
\end{align}
Since the matrix $A(y)$ is symmetric positive semidefinite for any $y\in\Rr^N$, we have 
\begin{equation}
\label{BEnorm}  
\|(I+\Delta t A(y))^{-1}\|=\dfrac{1}{1+\Delta t\lambda_{\min}(y)}\leqs 1,
\quad \forall y\in\Rr^N,
\end{equation}
with $\lambda_{\min}(y)\geqs 0$ being the smallest eigenvalue of $A(y)$.
Hence, by rewriting the backward Euler scheme as
$y^{n+1} = (I+\Delta t A(y^{n+1}))^{-1}(y^n + \Delta t g^{n+1})$ and taking the norm in the last
relation, we see that
the scheme~\eqref{BE} yields for any $\Delta t>0$ a bounded solution, i.e.,
\begin{equation}
\label{BEstab}
\|y^{n+1}\|\leqs\|y^n\| + \Delta t\|g^{n+1}\|,\quad n=0,1,2,\dots.   
\end{equation}

\begin{remark}
The boundedness properties stated in the paper imply 
stability for linear problems.
For discussion on stability for nonlinear
problems see, e.g., \cite[Sections~I.2.3, I.2.8]{HundsdorferVerwer:book}
and references therein.   
\end{remark}

Relation~\eqref{BE} is a system of nonlinear equations in~$y^{n+1}$.
To solve the system, in~\cite[App.~1, Ch.~2.11]{TikhonovSamarskii_UMF} the following iterative 
scheme is proposed:
\begin{equation}
\label{BEiter}
( I+\Delta t A(y^{(m)}) )y^{(m+1)} = y^n + \Delta t g^{n+1},
\quad m=0,1,2,\dots,   
\end{equation}
where the superscript $\cdot^{(m)}$ denotes the iteration number and we usually
take $y^{(0)}:=y^n$.  As argued in~\cite[Appendix~1]{TikhonovSamarskii_UMF},
the scheme is monotone for a certain finite-difference approximation
of one-dimensional heat equation~\eqref{heat}.  Below we
prove the monotonicity and convergence of the scheme for general, not necessarily
one-dimensional heat equation.  

\begin{corollary}
Assume the backward Euler method~\eqref{BE} in combination with iterative scheme~\eqref{BEiter}
is applied to solve the IVP~\eqref{ivp} with initial guess $y^{(0)}=y^n$ 
and relations~\eqref{aij},\eqref{g>=0} hold.
Then for all time steps $n=0,1,\dots$ iterative scheme~\eqref{BE},\eqref{BEiter}
\begin{enumerate}
\item 
is monotone, i.e., for any time step size $\Delta t>0$ and all iterations $m=0,1,\dots$   
\begin{equation}
\label{BEpos}
y^{(m+1)}\geqs 0,  
\end{equation}
\item
produces a sequence $y^{(m)}$ converging to solution $y^{n+1}$ of~\eqref{BE},
i.e., $\|y^{n+1}-y^{(m)}\|\rightarrow 0$ as $m\rightarrow\infty$,
provided that the time step size satisfies
\begin{equation}
\label{BEconv}
0<\Delta t < \frac{1}{L(\|y^n\| + \Delta t\|g^{n+1}\|)},      
\end{equation}
\item 
yields a bounded solution, i.e., for any time step size $\Delta t>0$ 
and all iterations $m=0,1,\dots$   
\begin{equation}
\label{BEst2}
\|y^{(m+1)}\|\leqs\|y^n\| + \Delta t\|g^{n+1}\|.
\end{equation}
\end{enumerate}
\end{corollary}

\begin{proof}
To prove monotonicity property~\eqref{BEpos}, we note that a matrix 
with nonpositive off-diagonal entries is positive
semidefinite if and only if 
it is a (possibly singular) $M$-matrix~\cite[Ch.~2.5]{HornJohnsonII}.  Therefore,
for any $\Delta t>0$ and any $y\in\Rr^N$ the matrix $I+\Delta t A(y)$ is a (nonsingular) 
$M$-matrix and its inverse is elementwise nonnegative: $(I+\Delta t A(y))^{-1}\geqs 0$.

To prove convergence of iterations~\eqref{BEiter}, we subtract relation~\eqref{BEiter}
from~\eqref{BE} and define the error vector $e^{(m)}\equiv y^{n+1}-y^{(m)}$.  Then
we obtain
\begin{gather*}
(I+\Delta t A(y^{n+1}))y^{n+1} - (I+\Delta t A(y^{(m)}))y^{(m+1)} = 0,
\\
e^{(m+1)} - \Delta t A(y^{(m)}) y^{(m+1)} = - \Delta t A(y^{n+1})y^{n+1},
\\
\intertext{and, adding $\Delta t A(y^{(m)}) y^{n+1}$ to both sides of the last equation,}
e^{(m+1)} + \Delta t A(y^{(m)}) y^{n+1} - \Delta t A(y^{(m)}) y^{(m+1)} = \Delta t A(y^{(m)}) y^{n+1} - \Delta t A(y^{n+1})y^{n+1},
\\
(I + \Delta t A(y^{(m)}))e^{(m+1)} = \Delta t ( A(y^{(m)}) - A(y^{n+1}))y^{n+1}.
\end{gather*}
Using~\eqref{BEnorm}, \eqref{L}, \eqref{BEstab}, we can bound
\begin{align*}
\|e^{(m+1)}\| &\leqs \Delta t \|( A(y^{(m)}) - A(y^{n+1}))y^{n+1}\| 
\leqs \Delta t \| A(y^{(m)}) - A(y^{n+1})\| \, \|y^{n+1}\| 
\leqs \Delta t L\| e^{(m)}\|\, \|y^{n+1}\| 
\\
&\leqs \Delta t L\| e^{(m)}\|\, \|y^{n+1}\| 
\leqs \Delta t L\| e^{(m)}\| (\|y^n\| + \Delta t\|g^{n+1}\|),
\end{align*}
from which convergence can easily be seen provided~\eqref{BEconv} holds.

Finally, the norm estimate~\eqref{BEst2} can be obtained in the same way as the estimate~\eqref{BEstab},
by taking the norm in the relations~\eqref{BEiter},\eqref{BEnorm}.
\end{proof}

In practice, the iterations~\eqref{BEiter} are stopped as soon as
$m>0$ (at least one iteration is done) and  
the residual of $y^{(m)}$ with respect to the nonlinear backward Euler 
equation~\eqref{BE}
$$
r^{(m)} = y^n + \Delta t g^{n+1} - ( I+\Delta t A(y^{(m)}) )y^{(m)}
$$ 
satisfies 
\begin{equation}
\label{BEstop}
\|r^{(m)}\| \leqs \texttt{tol} \cdot \|y^n + \Delta t g^{n+1}\|,
\end{equation}
with $\texttt{tol}>0$ being a given tolerance value.
We emphasize that the matrix in the residual expression is evaluated at $y^{(m)}$,  

\subsection{Nonlinear exponential Euler scheme}
In this method, a numerical solution $y^{n+1}\approx y(t_{n+1})$ approximating
solution $y(t)$ of~\eqref{ivp} at $t_{n+1}=(n+1)\Delta t$, $n=0,1,2,\dots$, is computed 
as follows.
\begin{equation}
\label{EE}
\begin{aligned}
&\text{solve IVP} \quad
\left\{\begin{matrix}
\tilde{y}'(t) = -A(\tilde{y}(t))\tilde{y}(t) + g^{n+1},
\quad t\in[t_n,t_n+\Delta t],
\\  
\tilde{y}(t_n) = y^n,
\end{matrix}\right.
\\
& \text{set}\quad y^{n+1}:=\tilde{y}(t_n+\Delta t),
\end{aligned}
\end{equation}
where $g^{n+1}=g(t_{n+1})$ and $y^0=v$.
We solve IVP in~\eqref{EE} iteratively,
by setting $\tilde{y}^{(0)}(t)\equiv y^n$, $t\in[t_n,t_n+\Delta t]$, and 
computing $\tilde{y}^{(m)}(t)\rightarrow\tilde{y}(t)$, $m=0,1,2,\dots$, $t\in[t_n,t_n+\Delta t]$. 
At each iteration $m=0,1,\dots$ we solve an IVP
\begin{equation}
\label{EEit}
\begin{aligned}
& (\tilde{y}^{(m+1)}(t))' = -\tilde{A}_m(t) \tilde{y}^{(m+1)}(t) + g^{n+1},
\quad t\in[t_n,t_n+\Delta t],
\\
& \tilde{y}^{(m+1)}(t_n) = y^n,
\end{aligned}
\end{equation}
where $\tilde{A}_m(t)=A(\tilde{y}^{(m)}(t))$.
Note that at the first iteration $m=0$ the matrix is constant since the initial
guess function $\tilde{y}^{(0)}(t)$ does not depend on time ($\tilde{y}^{(0)}(t)\equiv y^n$).
Method~\eqref{EE} can be seen as a nonlinear variant
of the exponential Euler scheme, see~\cite[relation~(1.6)]{HochbruckOstermann2010}.

We control convergence of the iterations~\eqref{EEit} by checking
residual 
of $\tilde{y}^{(m+1)}(t)$ with respect to ODE (ordinary differential equation) system
$\tilde{y}'(t) = -A(\tilde{y}(t))\tilde{y}(t) + g^{n+1}$:
\begin{align*}
  \tilde{r}^{(m+1)}(t) &= -A(\tilde{y}^{(m+1)}(t)) \tilde{y}^{(m+1)}(t) + g^{n+1} - (\tilde{y}^{(m+1)}(t))'
  \\
  &=
  -\tilde{A}_{m+1}(t) \tilde{y}^{(m+1)}(t) + g^{n+1}
  +\tilde{A}_m(t) \tilde{y}^{(m+1)}(t) - g^{n+1}
  = \left[ \tilde{A}_m(t) - \tilde{A}_{m+1}(t)\right] \tilde{y}^{(m+1)}(t).
\end{align*}
The iterations~\eqref{EEit} are stopped as soon as the residual is small in norm
at the final time $t=t_{n+1}$:
\begin{equation}
\label{EEstop}
\|\tilde{r}^{(m+1)}(t_{n+1})\| \leqs \texttt{tol} \cdot \|\tilde{A}_{m+1}(t_{n+1}) \tilde{y}^{(m+1)}(t_{n+1})\|,
\end{equation}
with $\texttt{tol}>0$ being a given tolerance value.

It turns out that in practice for typical time step sizes $\Delta t$, 
solution $\tilde{y}^{(m+1)}(t)$ of~\eqref{EEit} can be very well approximated 
by solution $y^{(m+1)}(t)$ of a simpler IVP
\begin{equation}
\label{EEit1}
\begin{aligned}
& (y^{(m+1)}(t))' = -A_m y^{(m+1)}(t) + g^{n+1},
\quad t\in[t_n,t_n+\Delta t],
\\
& y^{(m+1)}(t_n) = y^n,
\end{aligned}
\end{equation}
where the matrix is evaluated at $t=t_{n+1}$ and is kept constant: $A_m=A(y^{(m)}(t_{n+1}))$.
As the next corollary shows, the deviation $\|\tilde{y}^{(m+1)}(t) - y^{(m+1)}(t)\|$ between the
solutions of~\eqref{EEit} and~\eqref{EEit1} can be estimated and easily computed in practice.
To prove this result we first define an entire function
\begin{equation}
\label{phi}  
\varphi(z) \equiv \frac{e^z - 1}{z}, \quad z\in\Cc,
\end{equation}
with $\varphi(0)=1$, and formulate the following lemma.

\begin{lemma}
\cite[Section ~I.2.3]{HundsdorferVerwer:book}
Consider linear ODE system with variable
coefficients
$$
y'(t) = -A(t) y(t) + g(t) 
$$
and assume that there is a constant $\omega\in\Rr$ such that
for the matrix exponential $\exp(-sA)$ holds 
$\|\exp(-sA)\|\leqs e^{-s\omega}$, $s\in[0,T]$.  Then
\begin{equation}
\label{stab}
\begin{aligned}
\|y(t)\| 
&\leqs e^{-t\omega}\|y(0)\| + \int_0^t e^{-(t-s)\omega}\| g(s)\| ds \\
&\leqs e^{-t\omega}\|y(0)\| + t\varphi(-t\omega)\max_{s\in[0,t]}\| g(s)\|,
\quad t\in[0,T].
\end{aligned}
\end{equation}
\end{lemma}

\begin{proof}
See the last relation in~\cite[Section ~I.2.3]{HundsdorferVerwer:book}.
\end{proof}

Note that, for $t\geqs 0$, 
\begin{equation}
\label{phi_est}
t\varphi(-t\omega) =
\begin{cases}
\quad t, \quad &\omega = 0,
\\
\dfrac{1-e^{-t\omega}}{\omega},\quad &\omega>0.
\end{cases}  
\end{equation}

We now prove the result concerning solutions of~\eqref{EEit} and~\eqref{EEit1}.
\begin{corollary}
For all time steps $n=0,1,\dots$ and all iteration numbers $m=0,1,\dots$, 
solution $\tilde{y}^{(m+1)}(t)$ of~\eqref{EEit}
and solution $y^{(m+1)}(t)$ of~\eqref{EEit1} satisfy 
\begin{equation}
\label{errAconst}
\begin{aligned}
\|\tilde{y}^{(m+1)}(t_n+\tau) - y^{(m+1)}(t_n+\tau)\|\leqs
\tau\varphi(-\tau\omega) \max_{s\in[t_n,t_n+\tau]} 
\left\|\left[A_m-\tilde{A}_m(s)\right]y^{(m+1)}(s)\right\|,
\\
\tau\in[0,\Delta t],
\end{aligned}
\end{equation}
where $\omega$ is such that 
\begin{equation}
\label{omega}
\min_{(x,x)=1, x\in\Rr^N} (\tilde{A}_m(t)x,x)\geqs \omega \geqs 0, 
\quad \forall t\in[t_n,t_n+\Delta t].
\end{equation}
Note that $\omega\geqs 0$ due to the assumption that $A(y)$ is positive 
semidefinite for all $y\in\Rr^N$.
\end{corollary}

\begin{proof}
Note that for $m=0$ the matrix $\tilde{A}_m(t)$ is constant and equals $A_m$,
so that IVPs~\eqref{EEit} and~\eqref{EEit1} coincide and 
the estimate~\eqref{errAconst} trivially holds.
Consider the case $m>0$.
Without loss of generality we give a proof for the first time step $n=0$,
$t\in[0,\Delta t]$.
Substituting $y^{(m+1)}(t)$ in the ODE system 
$(\tilde{y}^{(m+1)}(t))' = -\tilde{A}_m(t) \tilde{y}^{(m+1)}(t) + g^{n+1}$,
we obtain a residual $\tilde{r}(t)$ of $y^{(m+1)}(t)$,
$$
\tilde{r}(t) = -\tilde{A}_m(t) y^{(m+1)}(t) + g^{n+1} - (y^{(m+1)}(t))' = 
\left[ A_m -\tilde{A}_m(t)\right] y^{(m+1)}(t),
$$
and see that $y^{(m+1)}(t)$ solves a perturbed ODE system
$$
(y^{(m+1)}(t))' = -\tilde{A}_m(t) y^{(m+1)}(t) + g^{n+1} - \tilde{r}(t),
\quad t\in[0,\Delta t].
$$
Subtracting this equation from the ODE system in~\eqref{EEit}, we arrive
at an IVP for error function $\tilde{e}(t)=\tilde{y}^{(m+1)}(t) - y^{(m+1)}(t)$:
\begin{equation}
\label{err_ivp}  
\tilde{e}'(t) = -\tilde{A}_m(t)\tilde{e}(t) + \tilde{r}(t), \quad \tilde{e}(t)=0,
\quad t\in[0,\Delta t].
\end{equation}
Condition~\eqref{omega} is equivalent to~\cite[Section~1.5]{Dekker-Verwer:1984},
\cite[Section~I.2.3]{HundsdorferVerwer:book}
\begin{equation}
\label{expA}
\|\exp(-t\tilde{A}_m(t))\| \leqs e^{-t\omega}, \quad t\geqs 0.
\end{equation}
Hence, applying the estimate~\eqref{stab} to~\eqref{err_ivp}, we have
$$
\|\tilde{e}(\tau)\| \leqs \tau\varphi(-\tau\omega)\max_{s\in[0,\tau]}
\left\|\left[ A_m -\tilde{A}_m(s)\right] y^{(m+1)}(s)\right\|, \quad \tau\in[0,\Delta t],
$$
which is the sought after estimate~\eqref{errAconst} for the first
time step $[0,\Delta t]$.
\end{proof}

Note that the error norm in~\eqref{errAconst} is zero at $t=t_n$ and the residual 
$\tilde{r}(t)$ is zero at $t=t_{n+1}$.  Therefore, to approximately evaluate 
the right hand side of the estimate~\eqref{errAconst} in practice, we can compute
$(s-t_n)\| [ A_m -\tilde{A}_m(s) ] y^{(m+1)}(s)\|$ at some 
point $s$ between $t_n$ and $t_n+\Delta t$, say at $s=t_n+\Delta t/2$.
For practical values $\Delta t>0$ this term often turns out to be negligibly small.
If this is not the case we can decrease the time step size $\Delta t$ or solve~\eqref{EEit1} 
successfully on smaller time subintervals covering $[t_n,t_n+\Delta t]$,
so that $\tilde{A}_m(t)$ is kept constant on these smaller subintervals
and the deviation $\|\tilde{y}^{(m+1)}(t) - y^{(m+1)}(t)\|$ becomes smaller.
In our experience, neglecting the large values 
$(s-t_n)\| [ A_m -\tilde{A}_m(s) ] y^{(m+1)}(s)\|$ is not disastrous and leads 
to an increase in the number of nonlinear iterations $m$.
In the following we assume that
$\Delta t>0$ is chosen such that $\|\tilde{y}^{(m+1)}(t) - y^{(m+1)}(t)\|$ is sufficiently small.
 
Since both the matrix $A_m$ and the source function 
$g^{n+1}$ are constant in~\eqref{EEit1}, iterations~\eqref{EEit1}
can be written in an equivalent form
\begin{equation}
\label{phi_it}
y^{(m+1)}(t_n+\tau) = y^n + \tau\varphi(-\tau A_m )(g^{n+1}- A_m y^{n}),
\quad \tau\in[0,\Delta t].
\end{equation}
Here $\varphi(-\tau A)$ is the matrix function defined by the $\varphi$ function~\eqref{phi},
see, e.g.,~\cite{Gantmacher,Higham_bookFM}.
In practice, we solve the IVP in~\eqref{EEit1} by computing 
$y^{(m+1)}(t_n+\Delta t)$ in~\eqref{phi_it} with the restarted Krylov
subspace method presented in~\cite{BKT21}.

The following result establishes monotonicity and convergence of the
iterative exponential Euler method~\eqref{EE},\eqref{EEit}.

\begin{corollary}
Assume the exponential Euler method~\eqref{EE} in combination with iterative 
scheme~\eqref{EEit}
is applied to solve the IVP~\eqref{ivp} 
and relations~\eqref{aij},\eqref{g>=0},\eqref{omega} hold.
For all time steps $n=0,1,\dots$ the iterative scheme~\eqref{EE},\eqref{EEit}
\begin{enumerate}
\item 
is monotone, i.e., for any time step size $\Delta t>0$ and all iterations $m=0,1,\dots$   
\begin{equation}
\label{EEpos}
\tilde{y}^{(m+1)}(t)\geqs 0,\quad t\in[t_n,t_n+\Delta t],  
\end{equation}
\item
produces a sequence $\tilde{y}^{(m)}(t)$ converging to solution $\tilde{y}(t)$ of~\eqref{EE},
i.e., $\max_{s\in[t_n,t_n+\Delta t]}\|\tilde{y}(s)-y^{(m)}(s)\|\rightarrow 0$ as 
$m\rightarrow\infty$, provided that the time step size $\Delta t$ satisfies
\begin{equation}
\label{EEconv}
0<\Delta t \varphi(-\Delta t\omega) L \max_{s\in[t_n,t_n+\Delta t]}\|\tilde{y}(s)\| < 1,
\end{equation}
where the function $\varphi$ is defined in~\eqref{phi},\eqref{phi_est} and 
$\omega$ in~\eqref{omega},
\item
yields a bounded solution, i.e., for any time step size $\Delta t>0$ 
and all iterations $m=0,1,\dots$   
\begin{equation}
\label{EEstab}
\|y^{(m+1)}(t_n+\tau)\|\leqs e^{-\tau\omega}\|y^n\| + \tau\varphi(-\tau\omega)\|g^{n+1}\|,
\quad \tau\in[0,\Delta t],
\end{equation}
where $\omega$ is defined in~\eqref{omega}.
\end{enumerate}  
\end{corollary}

\begin{figure}
\includegraphics[width=0.48\linewidth]{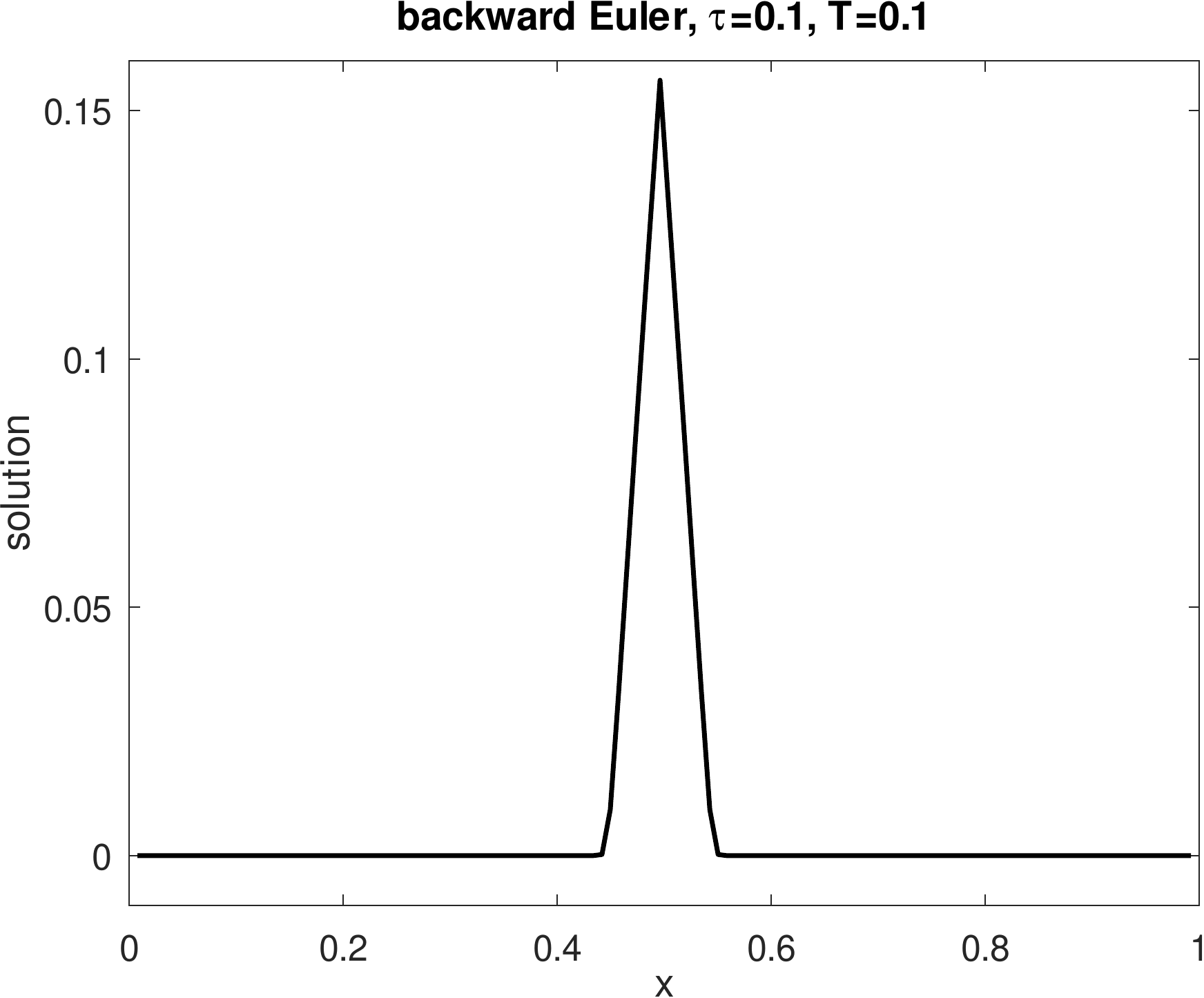}
\includegraphics[width=0.48\linewidth]{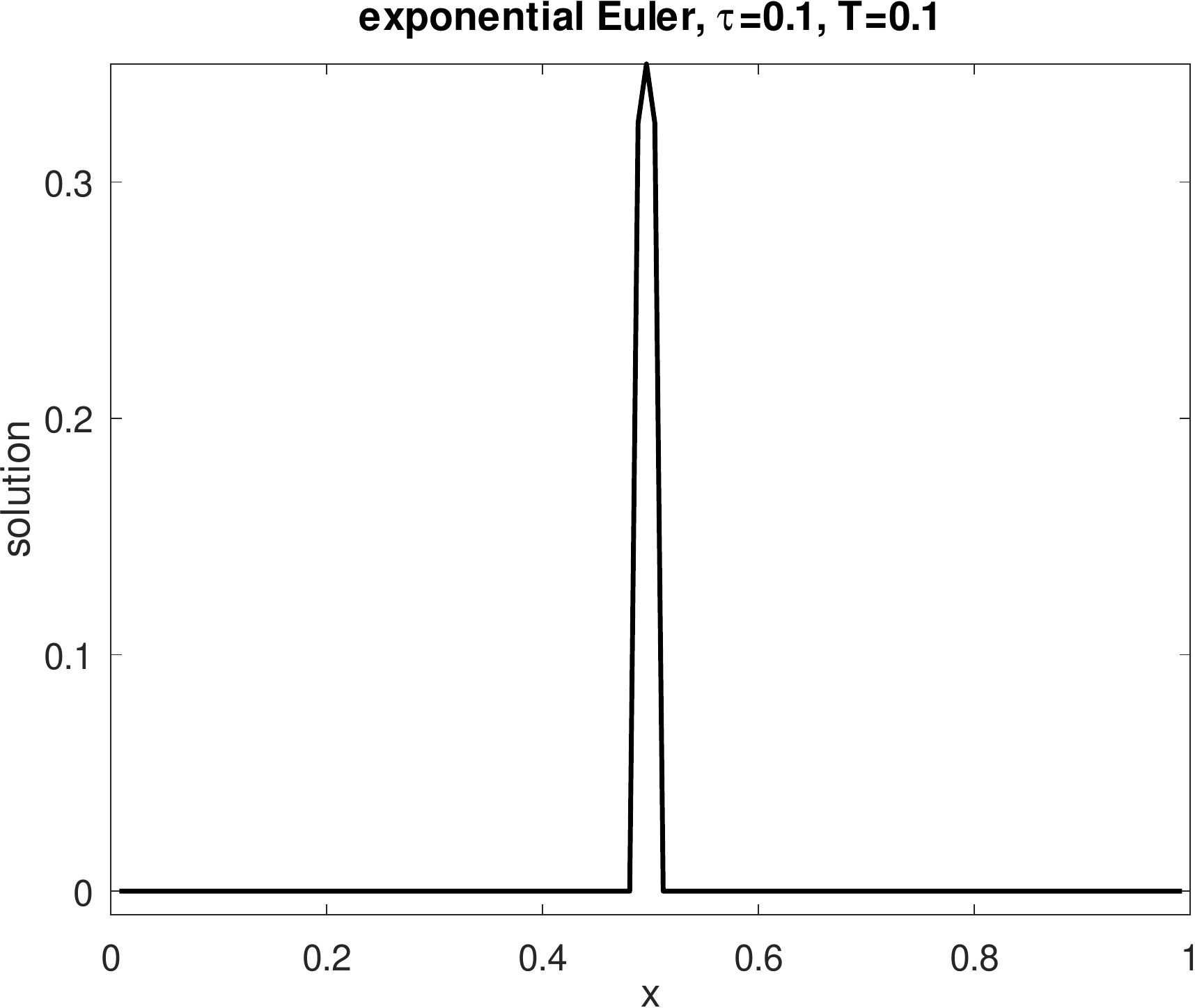}
\\[1ex]
\includegraphics[width=0.48\linewidth]{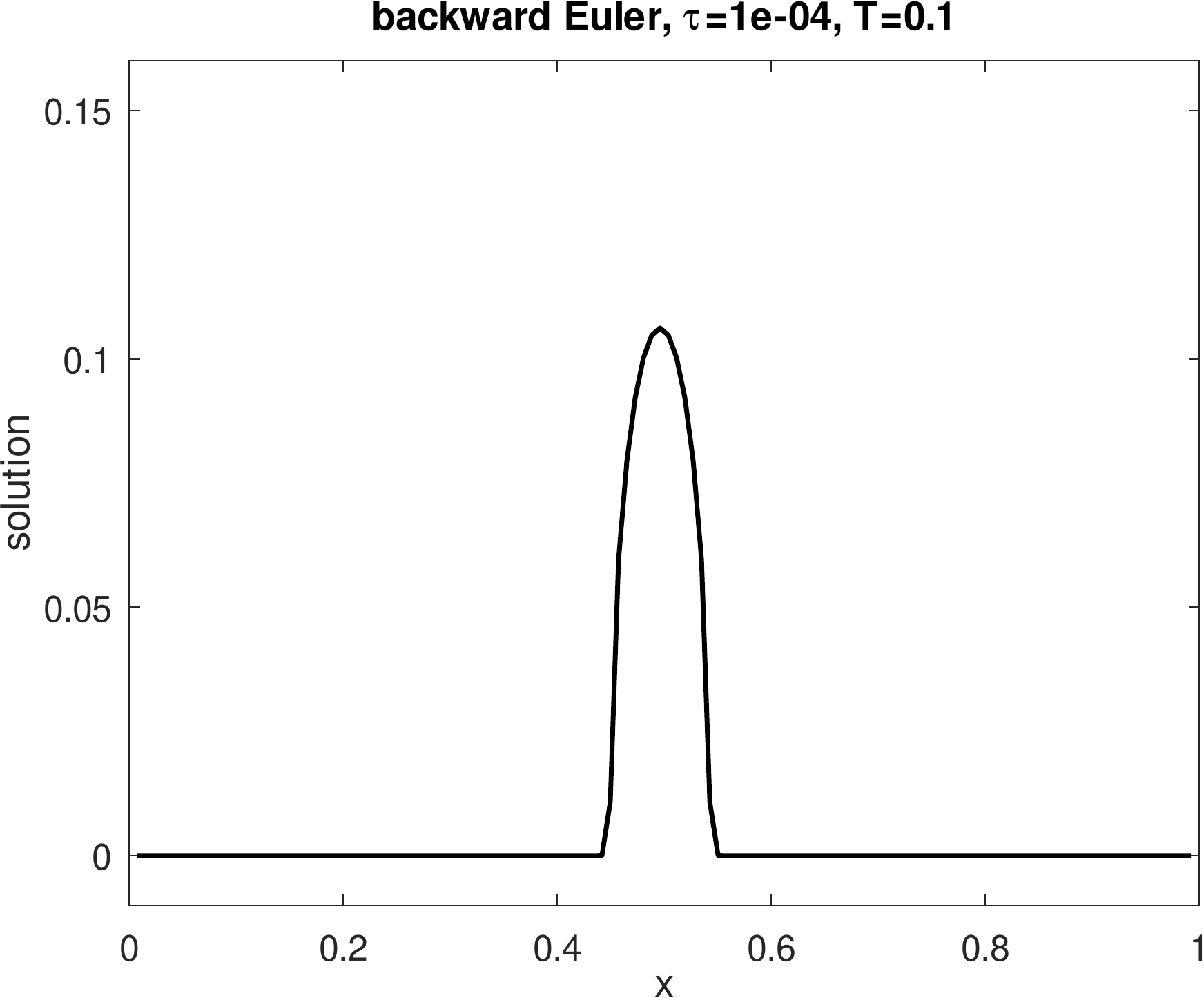}
\includegraphics[width=0.48\linewidth]{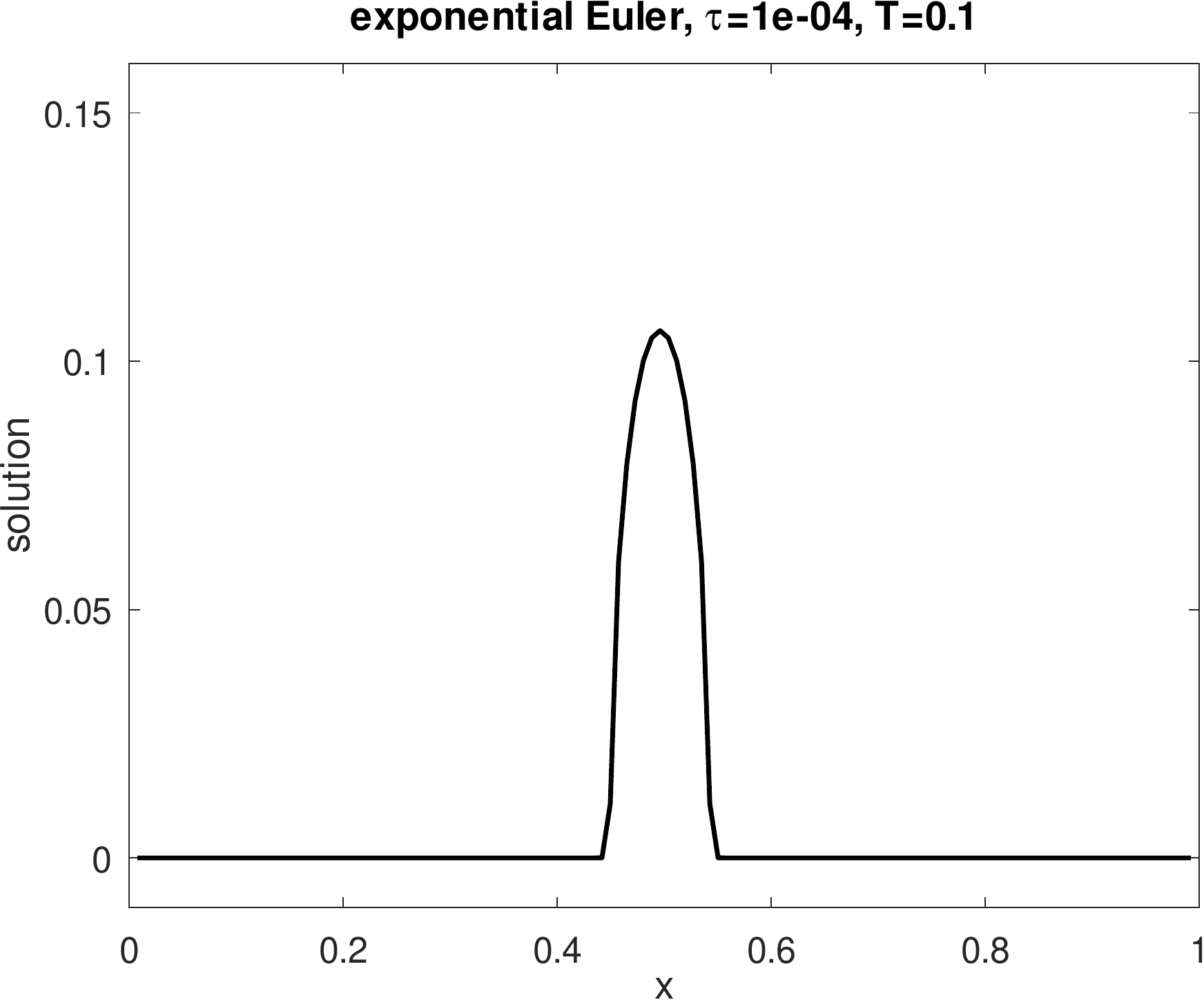}  
\caption{Numerical Green functions of the backward Euler (left) and exponential Euler (right) methods
at time $T=0.1$ computed with time steps $\Delta t=T$ (top) and $\Delta t=T/1000$ (bottom)
on a uniform grid with $N=128$ nodes.}
\label{f:Green}  
\end{figure}

\begin{proof}
Without loss of generality we give the proof for $n=0$, i.e.,
we consider the first time interval $[0,\Delta t]$ ($t_n=0$). 
However, to clearly see the connection between the derivations in the proof and 
the corollary, we keep on writing superindices containing $n$. 

The monotonicity of iterations~\eqref{EEit} can be established
by observing that the iterative approximations $\tilde{y}^{(m)}(t)$
solve IVP~\eqref{EEit}.  For this IVP we can show nonnegativity of the solution
in the same way as is done in the proof of Corollary~\ref{IVPpos}.

To prove convergence of the iterations, we subtract the ODE 
system~\eqref{EEit} from the ODE system~\eqref{EE} and obtain,
for error function $\epsilon^{(m+1)}(t)\equiv \tilde{y}(t) - y^{(m+1)}(t)$,
\begin{align*}
(\epsilon^{(m+1)}(t))' &= -A(\tilde{y}(t))\tilde{y}(t) + \tilde{A}_m(t) y^{(m+1)}(t), \\
(\epsilon^{(m+1)}(t))' &= -A(\tilde{y}(t))\tilde{y}(t) + \tilde{A}_m(t) y^{(m+1)}(t) 
                        -\tilde{A}_m(t)\tilde{y}(t) +\tilde{A}_m(t)\tilde{y}(t), \\
(\epsilon^{(m+1)}(t))' &= -\tilde{A}_m(t)\epsilon^{(m+1)}(t) +
\left[\tilde{A}_m(t)-A(\tilde{y}(t))\right]\tilde{y}(t),
\qquad t\in[0,\Delta t]. \\
\end{align*}
Applying to this ODE system estimate~\eqref{stab} and taking into
account initial condition $\epsilon^{(m+1)}(0)=0$ and
relations \eqref{omega}, \eqref{expA}, we get
\begin{align}
\|\epsilon^{(m+1)}(t)\| 
&\leqs t\varphi(-t\omega) \, \max_{s\in[0,t]}
\left\|\left[ \tilde{A}_m(s)-A(\tilde{y}(s))\right]\tilde{y}(s)\right\|
\notag\\
&= t\varphi(-t\omega) \, \max_{s\in[0,t]}
\left\|\left[ A (y^{(m)}(s))-A(\tilde{y}(s))\right]\tilde{y}(s)\right\|
\notag\\\displaybreak[2]
&\leqs t\varphi(-t\omega) L \max_{s\in[0,t]}\bigl\|
\underbrace{y^{(m)}(s) -\tilde{y}(s)}_{\epsilon^{(m)}(s)}\bigr\|
\max_{s\in[0,t]}\|\tilde{y}(s)\|,  
\quad t\in[0,\Delta t], 
\end{align}
where we use definition of $\tilde{A}_m(t)$ and property~\eqref{L}.
Since $t\varphi(-t\omega)$ is a monotonically increasing function, we have
$$
\max_{s\in[0,\Delta t]}\|\epsilon^{(m+1)}(s)\|\leqs
\Delta t\varphi(-\Delta t\omega) L \max_{s\in[0,\Delta t]}\|\epsilon^{(m)}(s)\|
\max_{s\in[0,\Delta t]}\|\tilde{y}(s)\|,
$$ 
which implies $\max_{s\in[0,\Delta t]}\|\epsilon^{(m+1)}(s)\|<\max_{s\in[0,\Delta t]}\|\epsilon^{(m)}(s)\|$
provided~\eqref{EEconv} holds.

The norm bound~\eqref{EEstab} can be obtained by applying
estimate~\eqref{stab} to IVP~\eqref{EEit}.
\end{proof}

\begin{remark}
Taking into account~\eqref{phi_est},
we see that the convergence conditions~\eqref{BEconv} and~\eqref{EEconv} 
for the backward and exponential Euler schemes are very similar.  
The convergence condition for the exponential Euler scheme is 
less restrictive if $\omega > 0$.
\end{remark}

\section{Numerical experiments}
\subsection{1D heat equation}
This test is considered in~\cite[App.~1, Ch.~2.11]{TikhonovSamarskii_UMF}.
It is a one-dimensional problem~\eqref{heat} in domain $\Omega = [0,1]$
with $k_0 = 0.5$, $\sigma  = 2$ and exact solution
\begin{equation}
\label{uex1D}
u_{\mathrm{exact}}(x,t) = 
\begin{cases}
\left(\dfrac{\sigma c}{k_0}
(ct - x) \right)^{1/\sigma}, \quad &\text{for} \quad x\leqs ct,
\\
0,   \quad &\text{for} \quad x> ct,
\end{cases}
\end{equation}
where $c$ is the heat wave speed (set to $c=1$ in this test).
We take initial and boundary conditions in~\eqref{heat} 
consistent with the exact solution and solve the problem on
the time interval $0\leqs t\leqs T=0.5$.
The error values reported below are computed as
\begin{equation}
\label{err}
\frac{\|y^n - y_{\mathrm{exact}}(T)\|}{\|y_{\mathrm{exact}}(T)\|}, \quad 
n\Delta t = T,
\end{equation}
where the vector $y_{\mathrm{exact}}(T)$ contains the grid values
of the exact solution function at final time $T$.
Furthermore, we set the accuracy tolerance in~\eqref{BEstop},\eqref{EEstop} 
to $\texttt{tol}=10^{-2}$ 
and the tolerance for computing actions of the $\varphi$ matrix function
to $10\,\texttt{tol}$.  Taking more stringent tolerance values do not give
a higher accuracy, as the error is determined by spatial discretization.
The restarted Krylov subspace solver for evaluating the $\varphi$ matrix function
is employed with the maximal Krylov subspace dimension~$30$
(in most runs no restarts are necessary).

We first test monotonicity of both schemes by computing 
numerical solution for homogeneous Dirichlet boundary conditions
and initial vector $v$ with a single entry set to~1, 
and all the other entries being zero (cf.~\eqref{ivp}).
The solution to this problem, which can be seen as a numerical
equivalent of the Green function, is shown in Figure~\ref{f:Green}.
As we see, both backward Euler method and exponential Euler
method manage to produce nonnegative solutions for this test.  
 
\begin{figure}
\includegraphics[width=0.48\linewidth]{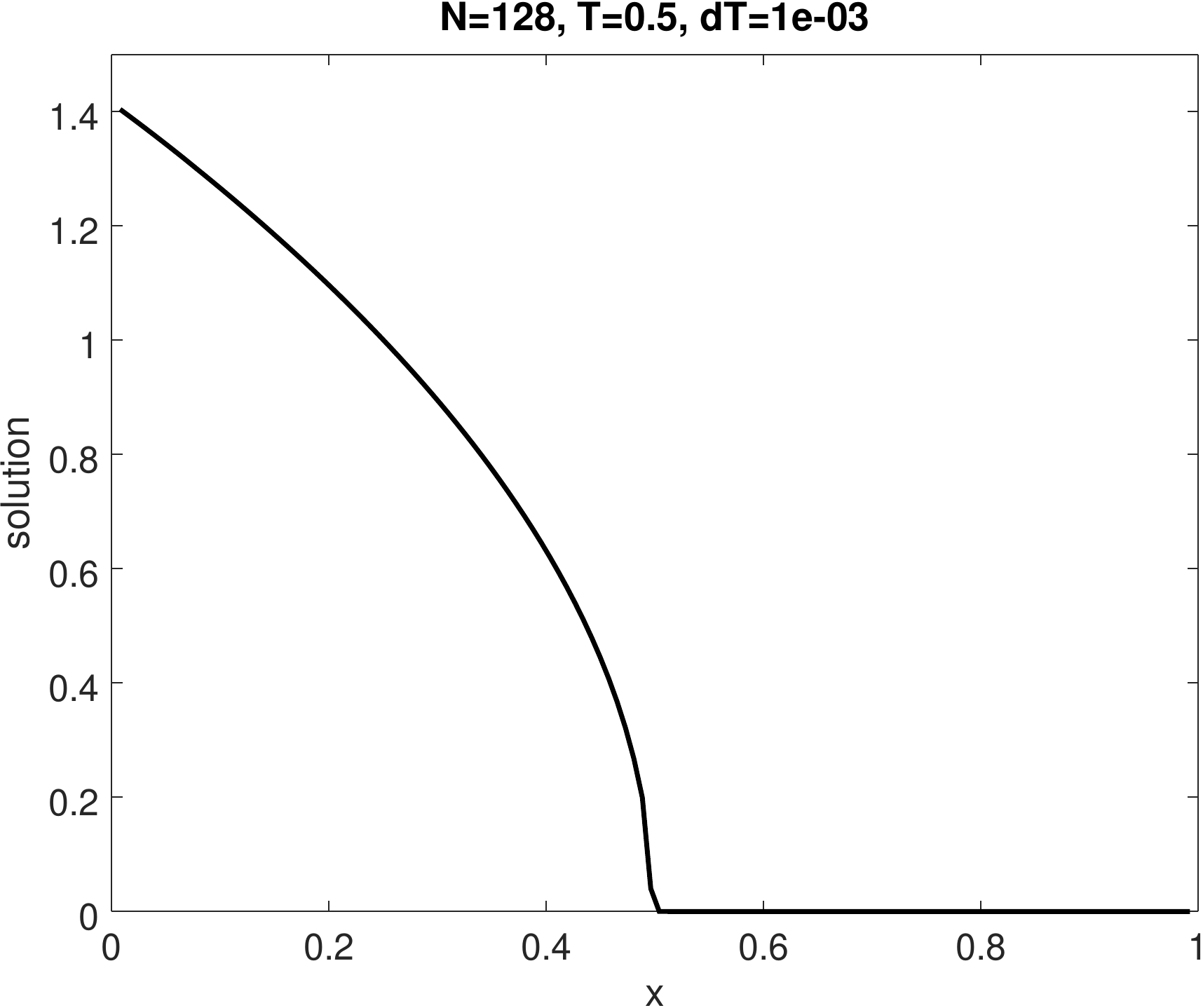}
\includegraphics[width=0.48\linewidth]{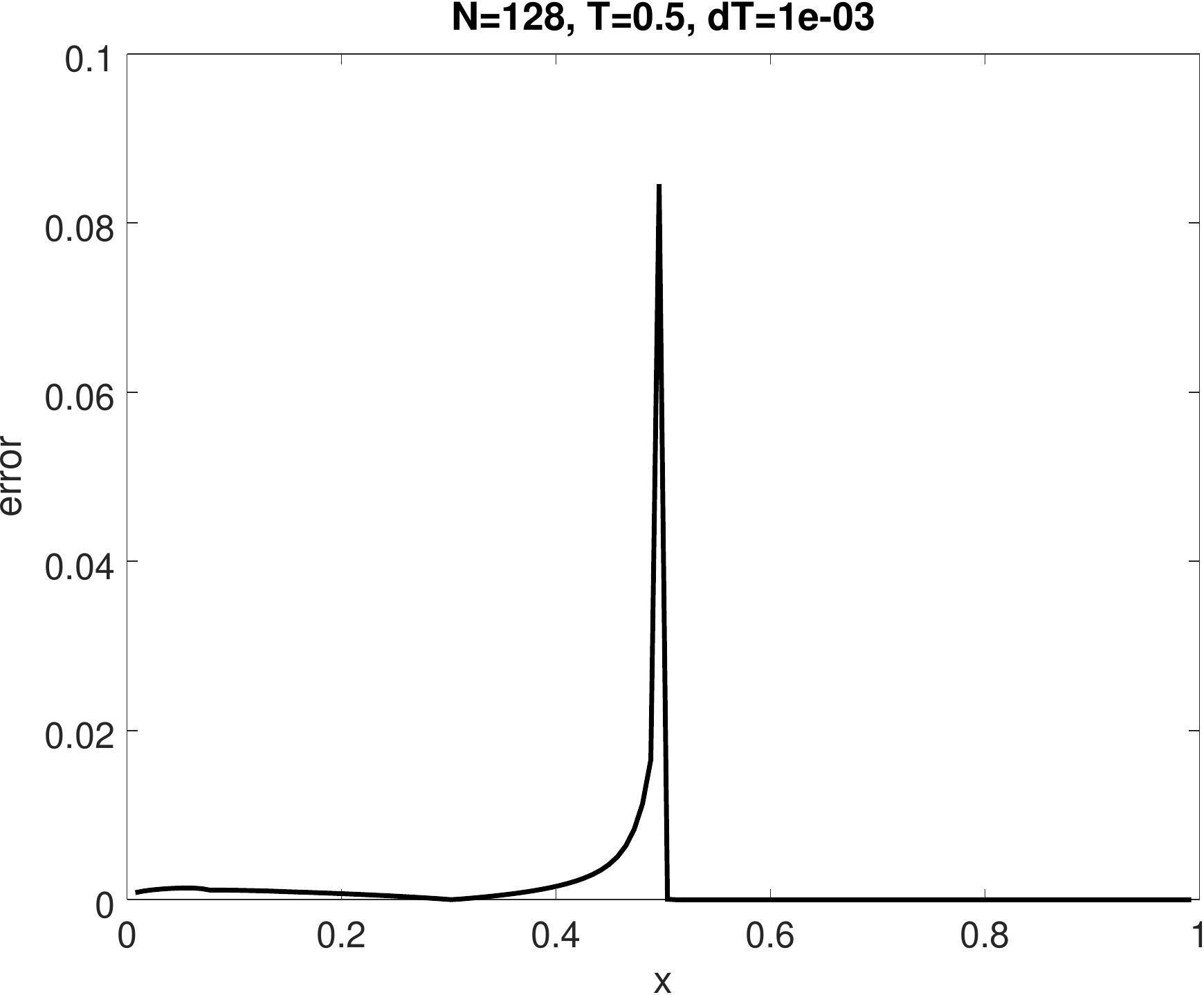}
\caption{Solution (left) and error (right) of the exponential Euler scheme on
the grid $N=128$ at $T=0.5$, computed with the time step size $\Delta t=10^{-3}$
($\Delta t \max_{t\in[0,T]}\|A(y(t))\|_1\approx 65$).}
\label{f:1D}  
\end{figure}

We now consider the regular test runs, with the boundary and initial conditions
determined by the exact solution~\eqref{uex1D}. 
A numerical solution and the corresponding error for the exponential Euler
method are shown in Figure~\ref{f:1D}.

\addtolength{\tabcolsep}{2ex}
\begin{table}
\caption{1D heat equation test problem. 
Results for the backward Euler (BE) and exponential Euler (EE) methods.  
The smallest matvec number for both methods and on each grid is underlined.}
\label{t:res1D}

\bigskip

\centerline{\begin{tabular}{ccc}
\hline\hline
$\Delta t$ & BE & EE \\
       & \texttt{error}, \#{} iterations & \texttt{error}, \#{} iterations\\
       &   (\#{} matvecs)                       &   (\#{} matvecs)      \\ 
\hline
\multicolumn{3}{c}{grid $N=128$, $\max_{t\in[0,T]}\|A(y(t))\|_1\approx 6.5\cdot 10^4$} \\
{\tt5e-05} & {\tt4.07e-03}, 10000 (35646)  & {\tt5.26e-03}, 10039 (10926)\\
{\tt1e-04} & {\tt4.49e-03},  5000 (20206)  & {\tt5.63e-03},  5072  (8053)\\
{\tt5e-04} & {\tt7.95e-03},  1008 (7098)   & {\tt9.08e-03},  1118  (4247)\\
{\tt1e-03} & {\tt1.18e-02},   587 (\underline{5102})   & {\tt1.11e-02},   642  (\underline{3473})\\
\hline
\multicolumn{3}{c}{grid $N = 256$, $\max_{t\in[0,T]}\|A(y(t))\|_1\approx 2.6\cdot 10^5$} \\
{\tt5e-05} & {\tt2.33e-03}, 10000 (50130)  & {\tt3.61e-03}, 10073 (20584)\\
{\tt1e-04} & {\tt2.84e-03}, 5000 (32482)   & {\tt4.66e-03}, 5102 (15241)\\
{\tt5e-04} & {\tt6.65e-03}, 1086 (13276)   & {\tt7.57e-03}, 1142  (9318)\\
{\tt1e-03} & {\tt1.12e-02},  668 (\underline{9980})    & {\tt1.08e-02},  768 (\underline{7526})\\
\hline
\end{tabular}}
\end{table}
\addtolength{\tabcolsep}{-2ex}

In Table~\ref{t:res1D} nonlinear iteration numbers are shown for both 
methods for various spatial grids and time step sizes.  
For each grid 
the value of $\max_{t\in[0,T]}\|A(y(t))\|_1$ is also shown in the table.
The matvec (matrix-vector product) numbers for the exponential Euler method
are the total numbers of the Krylov subspace iterations required
to evaluate the $\varphi$ matrix functions.
The matvec values reported in the table for the backward Euler method have 
the following meaning.
Solution of the linear system~\eqref{BEiter} in the backward Euler 
method can be replaced by a certain number of Chebyshev iterations.
The number of Chebyshev iterations can then be chosen such that the scheme
is stable (for linear problems) eventhough the linear systems~\eqref{BEiter}
are solved approximately. 
This leads to the so-called explicit local iteration (LI) schemes, which are 
known to work well for heat conduction problems~\cite{Zhukov2011}.
In the table the bracketed matvec values for the backward Euler method
are the numbers of Chebyshev iterations needed in the monotone LI-M
scheme~\cite{Zhukov2011}.  As Chebyshev iterations require a singe matvec per iteration,
the shown Chebyshev iteration numbers equal the matvec numbers in the LI-M method.
We see that the exponential Euler scheme is more efficient in terms
of the matvecs than the LI-M method.

\subsection{2D heat equation}
We solve problem~\eqref{heat} in domain $\Omega = [0,1]\times [0,1]$ and 
we take $k_0=1$ and $\sigma=2$ in~\eqref{k}.  The initial and boundary 
conditions in~\eqref{heat} are
taken such that problem~\eqref{heat} has exact solution
\begin{equation}
\label{uex}
u_{\mathrm{exact}}(x, y, t) = \dfrac{1}{\sqrt[3]{t}}
\sqrt{\dfrac{\max\left\{ 0; \left(1.3 - \frac{x^2 +y^2}{\sqrt[3]{t}}\right) \right\} }{6} }
\,.
\end{equation}
The time interval is chosen to be $[t_0,T]=[0.0001, 0.0051]$.
The error values reported below are computed
according to~\eqref{err}.
Just as in the previous test, in~\eqref{BEstop},\eqref{EEstop} 
the nonlinear accuracy tolerance is set to $\texttt{tol}=10^{-2}$ 
and the tolerance for computing actions of the $\varphi$ matrix function
to $10\,\texttt{tol}$.
The maximal Krylov subspace dimension for evaluating the $\varphi$ matrix function
is again~$30$.

\addtolength{\tabcolsep}{2ex}
\begin{table}
\caption{2D heat equation test problem.
Results for the backward Euler (BE) and exponential Euler (EE) methods.  
The smallest matvec number for both methods and on each grid is underlined.}
\label{t:res2D}

\bigskip

\centerline{\begin{tabular}{ccc}
\hline\hline
$\Delta t$ & BE & EE \\
       & \texttt{error}, \#{} iterations & \texttt{error}, \#{} iterations\\
       &   (\#{} matvecs)                       &   (\#{} matvecs)      \\ 
\hline
\multicolumn{3}{c}{grid $64\times 64$, $\max_{t\in[t_0,T]}\|A(y(t))\|_1\approx 3.0\cdot10^6$} \\
{\tt1e-06} & {\tt1.24e-02}, 5000 (12306)  & {\tt1.20e-02}, 5000 (5079)\\
{\tt5e-06} & {\tt1.18e-02}, 1011  (4266)  & {\tt1.16e-02}, 1038 (\underline{1601})\\
{\tt1e-05} & {\tt1.17e-02},  536  (\underline{2668})  
                                          & {\tt1.17e-02}, 613 (1806)\\
{\tt5e-05} & {\tt1.84e-02},  238 (2884)   & {\tt1.75e-02}, 479 (4164)\\
\hline
\multicolumn{3}{c}{grid $128\times 128$, $\max_{t\in[t_0,T]}\|A(y(t))\|_1\approx 1.3\cdot10^7$} \\
{\tt1e-06} & {\tt7.44e-03}, 5000 (20556)  & {\tt7.20e-03}, 5000 (6135)\\
{\tt5e-06} & {\tt7.24e-03}, 1032  (6654)  & {\tt7.51e-03}, 1103 (\underline{4129})\\
{\tt1e-05} & {\tt7.65e-03},  613  (\underline{5814})  
                                          & {\tt8.52e-03}, 732 (5210)\\
{\tt5e-05} & {\tt1.95e-02},  398  (9688)  & {\tt2.51e-02}, 1045 (18403)\\
\hline
\multicolumn{3}{c}{grid $256\times 256$, $\max_{t\in[t_0,T]}\|A(y(t))\|_1\approx 4.9\cdot10^7$} \\
{\tt1e-06} & {\tt3.13e-03}, 5000 (27700)  & {\tt3.34e-03}, 5000 (12746)\\
{\tt5e-06} & {\tt4.44e-03}, 1082 (\underline{12860})  
                                          & {\tt5.88e-03}, 1202 (\underline{10800})\\
{\tt1e-05} & {\tt6.75e-03},  699 (13728)  & {\tt9.51e-03}, 1067 (17955)\\
{\tt5e-05} & {\tt2.22e-02},  710 (34812)  & {\tt3.76e-02}, 2759 (106906)\\
\hline
\end{tabular}}
\end{table}
\addtolength{\tabcolsep}{-2ex}

In Table~\ref{t:res2D} the results for this 2D test are presented in the same
form as for the previous test in Table~\ref{t:res1D}.
The results show the similar trend: the exponential Euler method
appears to be more efficient in terms of required matvec numbers than
the LI-M method (i.e., the backward Euler scheme combined with Chebyshev
iterations to solve~\eqref{BEiter}).    

\section{Conclusions and an outlook to further research}
For nonlinear heat conduction problems
we have compared, theoretically and numerically, nonlinear iterative
methods based on the backward Euler and on the exponential Euler 
schemes.  Both methods are shown to produce monotone and bounded 
solutions and the underlying nonlinear iterations appear to have
similar convergence properties.  In the experiments 
both methods exhibit a robust performance.
The proposed nonlinear exponential Euler method is an explicit 
scheme.  Earlier, explicit schemes based on Chebyshev iterations
(the so-called local iteration schemes) 
have been shown to work successfully for nonlinear heat conduction 
problems~\cite{Zhukov2011}.
Preliminary estimates of computational work presented here suggest
that the nonlinear exponential Euler outperforms
the LI-M (local iteration modified) scheme.  However, actual comparisons
of the proposed exponential Euler method and the local
iterations methods have yet to be done. 

\begin{acknowledgments}
The work of the first author is supported by the Russian Science
Foundation grant~19-11-00338.
The first author thanks Leonid Knizhherman for a useful suggestion.
\end{acknowledgments}

\medskip

This is a preprint of the work accepted for publication in
Lobachevskii Journal of Mathematics (\url{https://www.pleiades.online/})
\copyright~2022, the copyright holder indicated in the
Journal.

\bibliography{matfun,my_bib}

\end{document}